\documentclass[12]{amsart}
\usepackage{amsmath,amssymb,amsthm,color,enumerate,comment,centernot,enumitem,url,cite}
\usepackage{graphicx,relsize,bm}
\usepackage{mathtools}
\usepackage{array}

\makeatletter
\newcommand{\tpmod}[1]{{\@displayfalse\pmod{#1}}}
\makeatother

\newtheorem{thm}{Theorem}[section]

\theoremstyle{remark}

\theoremstyle{definition}

\newtheorem{rem}[thm]{Remark}

\theoremstyle{THM}

\newcommand{\abs}[1]{\left|{#1}\right|}

\def\FF {{\mathcal F}}

\def\Z {{\mathbb Z}}

\def\Q {{\mathbb Q}}

\def\Z {{\mathbb Z}}
\def\Q {{\mathbb Q}}

\def\Gal{{\mbox{{\rm{Gal}}}}}

\makeatletter
\@namedef{subjclassname@2020}{%
  \textup{2020} Mathematics Subject Classification}
\makeatother

\def\red#1 {\textcolor{red}{#1 }}
\def\blue#1 {\textcolor{blue}{#1 }}

\numberwithin{equation}{section}

\def\Z {{\mathbb Z}}

\begin{document}

\title[Monogenic trinomials and class numbers of related quadratic fields]{Monogenic trinomials and class numbers of related quadratic fields}


\author{Lenny Jones}
\address{Professor Emeritus, Department of Mathematics, Shippensburg University, Shippensburg, Pennsylvania 17257, USA}
\email[Lenny~Jones]{doctorlennyjones@gmail.com}

\date{\today}

\begin{abstract}
We say that a monic polynomial $f(x)\in \Z[x]$ of degree $N\ge 2$ is \emph{monogenic} if $f(x)$ is irreducible over $\Q$ and $\{1,\theta,\theta^2,\ldots ,\theta^{N-1}\}$ is a basis for the ring of integers of $\Q(\theta)$, where $f(\theta)=0$. In this article, we investigate the divisibility of the class numbers of quadratic fields $\Q(\sqrt{\delta})$ for certain families of monogenic trinomials $f(x)=x^N+Ax+B$, where $\delta\ne \pm 1$ is a squarefree divisor of the discriminant of $f(x)$.
\end{abstract}


\maketitle

\section{Introduction}\label{Section:Intro}
We say that a monic polynomial $f(x)\in \Z[x]$ is \emph{monogenic} if $f(x)$ is irreducible over $\Q$ and $\{1,\theta,\theta^2,\ldots ,\theta^{\deg{f}-1}\}$ is a basis for the ring of integers of $\Q(\theta)$, where $f(\theta)=0$. Hence, $f(x)$ is monogenic if and only if  $\Z_K=\Z[\theta]$. For the minimal polynomial $f(x)$ of an algebraic integer $\theta$ over $\Q$, it is well known \cite{Cohen} that
\begin{equation} \label{Eq:Dis-Dis}
\Delta(f)=\left[\Z_K:\Z[\theta]\right]^2\Delta(K),
\end{equation}
where $\Delta(f)$ and $\Delta(K)$ are the respective discriminants over $\Q$ of $f(x)$ and the number field $K$.
Thus, from \eqref{Eq:Dis-Dis}, $f(x)$ is monogenic if and only if $\Delta(f)=\Delta(K)$.

Many authors have investigated divisibility properties of class numbers of quadratic fields. We present a few such results from the literature that are germane to our investigations in this paper.
The first theorem is from 1955 and is due to Ankeny and Chowla \cite{AC}.
\begin{thm}\label{Thm:AC}
  Let $n$ and $M$ be positive integers with $M\ge 5$. If $\delta:=M^{2n}+1$ is squarefree, then $n$ divides the class number $h_K$ of the real quadratic field $K=\Q(\sqrt{\delta})$.
\end{thm}
The next theorem is due to M. Ram Murty \cite{Murty}, which he published in 1998.
\begin{thm}\label{Thm:Murty}
Suppose that $\delta:=1-M^n$ is squarefree with $M\ge 5$ and odd. Then the class group of $\Q(\sqrt{\delta})$ has an element of order $n$.
\end{thm}
\begin{rem}\label{Rem:Murty}
  We point out that Murty's proof of Theorem \ref{Thm:Murty} is also valid when $M=3$ and $n$ is odd.
\end{rem}
In 2000, Kishi and Miyake \cite{KM} gave the following parameterization of the quadratic fields whose class numbers are divisible by 3:
 \begin{thm}\label{Thm:KM}
 Let $u,w\in \Z$ with $\gcd(u,w)=1$, such that $d:=4uw^3-27u^2$ is not a square and one of the following sets of conditions holds:
 \begin{align}
\label{C1}  &3\nmid w;\\
\label{C2}  &3\mid w, \quad uw\not \equiv 3 \pmod{9}, \quad u\equiv w \pm 1 \pmod{9};\\
\label{C3}  &3\mid w, \quad uw \equiv 3 \pmod{9}, \quad u\equiv w\pm 1 \pmod{27}.
 \end{align}
 Define
 \begin{equation}\label{Eq:f}
 f_{u,w}(x):=x^3-uwx-u^2.
 \end{equation} If $f_{u,w}(x)$ is irreducible over $\Q$, then the roots of $f_{u,w}(x)=0$ generate an unramified cyclic cubic extension of the quadratic field $\Q(\sqrt{d})$. Conversely, every quadratic field $K$ whose class number is divisible by {\rm 3} and every unramified cyclic cubic extension of $K$ are given in this way by a suitable pair of integers $u$ and $w$.
 \end{thm}

The next theorem is from 2009 and due to Kishi \cite{K}.
\begin{thm}\label{Thm:K}
For any positive integers $k$ and $n$ with $2^{2k}<3^n$ and $(k,n)\ne (2,3)$, the class number of the imaginary quadratic field $\Q(\sqrt{2^{2k}-3^n})$ is divisible by $n$.
\end{thm}

In this paper, for certain monogenic trinomials $f(x)=x^N+Ax+B$, we use Theorems \ref{Thm:Murty}, \ref{Thm:KM} and \ref{Thm:K} to investigate divisibility properties of the class numbers of quadratic fields $\Q(\sqrt{\delta})$, where $\delta\ne \pm 1$ is a divisor of $\Delta(f)$. One motivation for such investigations is to determine whether there is some possible connection between the monogenicity of $f(x)$ and the class number of $\Q(\sqrt{\delta})$. This notion is not completely outside the realm of possibility since there are well-known connections between $\Delta(f)$ and the monogenicity of $f(x)$. For example, if $f(x)$ is monic and irreducible over $\Q$ such that $\Delta(f)$ is squarefree, then $f(x)$ is monogenic. A second motivation is to prove that there exist infinitely many distinct monogenic polynomials such that the class number of the related quadratic field can be arbitrarily large.

We summarize our approach to this study. Our somewhat natural starting point is to construct and examine families of polynomials $f(x)$ such that $\Delta(f)$ has a divisor that is of the same form as the discriminants for the quadratic fields found in one of the known results: Theorem \ref{Thm:Murty}, Theorem \ref{Thm:KM} or Theorem \ref{Thm:K}. For the construction of such polynomials, we focus on trinomials $f(x)=x^N+Ax+B$ because of the relative simplicity of their discriminants. Then, for each of these constructed families of trinomials, we provide necessary and sufficient conditions for the monogenicity of a trinomial in the family, and we show that no two such degree-$N$ monogenic trinomials generate the same extension over $\Q$.  Next, we confirm that these monogenic trinomials do satisfy the conditions of one of the Theorems \ref{Thm:Murty}, \ref{Thm:KM} and \ref{Thm:K}, so that the divisibility conditions for the related quadratic field found in these known theorems do indeed hold. Finally, we examine the data for any obvious patterns.

Although our investigations did not yield any apparent connections between the monogenicity of these trinomials and the divisibility of the class numbers of the related quadratic fields, we have been able to show that there exist infinitely many monogenic trinomials that fall within the context of each of the Theorems \ref{Thm:Murty}, \ref{Thm:KM} and \ref{Thm:K}. Moreover, as mentioned above, we also show that no two degree-$N$ monogenic trinomials generate the same extension over $\Q$. Consequently, in doing so, we have established our second motivation. 

In the statement of each new theorem below, we indicate explicitly its association to the particular one of Theorems \ref{Thm:Murty}, \ref{Thm:KM} and \ref{Thm:K}.
 More precisely, our main results are:
\begin{thm}(associated with Theorem \ref{Thm:KM})\label{Thm:Main1}\\
Let $u$, $w$ and $d$ be integers that satisfy the conditions set forth in Theorem \ref{Thm:KM}, and let $f_{u,w}(x)$ be as defined in \eqref{Eq:f}. Suppose that $f_{u,w}(x)$ is irreducible over $\Q$. Then
\begin{equation}\label{Eq:Mono1}
f_{u,w}(x)\ \mbox{is monogenic if and only if $u=1$ and $d$ is squarefree,}
\end{equation} the set
\begin{equation}\label{FF1}
\FF_1=\{f_{1,w}(x): \mbox{$d$ is squarefree},\}
\end{equation} is infinite, every element of $\FF_1$
is such that $\Gal_{\Q}(f_{a,b})\simeq S_3$, the class number of the unique quadratic subfield $\Q(\sqrt{d})$ of the splitting field of $f_{a,b}(x)$ is divisible by $3$,  and no two elements of $\FF_1$ generate the same cubic extension over $\Q$.
 \end{thm}
 \begin{thm}\label{Thm:Main2}(associated with Theorem \ref{Thm:KM})
   Let $a$ and $b$ be positive integers with $b\ge 2$ and $(a,b)\ne (1,2)$. Define
   \[f_{a,b}(x):=x^3-2^{2a}x-3^{b-1} \quad \mbox{and} \quad \delta:=2^{6a+2}-3^{2b+1}.\]
   Then
   \begin{equation}\label{Eq:Mono2}
f_{a,b}(x)\ \mbox{is monogenic if and only if $\delta$ is squarefree,}
\end{equation} and every element of the set
\begin{equation}\label{FF2}
\FF_2=\{f_{a,b}(x): \mbox{$\delta<0$ is squarefree}\}
\end{equation} is such that $\Gal_{\Q}(f_{a,b})\simeq S_3$, $\FF_1\cap \FF_2=\varnothing$, no two elements of $\FF_2$ generate the same cubic extension of $\Q$, and the class number of the unique imaginary quadratic subfield $\Q(\sqrt{\delta})$ of the splitting field of $f_{a,b}(x)$ is divisible by $(6b+3)/\gcd(3,2b+1)$.
 \end{thm}

\begin{thm}\label{Thm:Main3}(associated with Theorem \ref{Thm:AC})
Let $b,N\in \Z$ such that $b\ge 1$, $N\ge 6$ and $N\equiv 2 \pmod{4}$.
Define
 \[f_{N,b}(x):=x^N-N(N-1)^bx-(N-1) \quad \mbox{and} \quad \delta:=(N-1)^{bN}+1.\]
 Suppose that $f_{N,b}(x)$ is irreducible over $\Q$. Then
\begin{equation}\label{Eq:Mono3}
f_{N,b}(x)\ \mbox{is monogenic if and only if $N(N-1)$ and $\delta$ are squarefree.}
\end{equation} Furthermore, no two such monogenic trinomials $f_{N,b}(x)$ generate the same $N$th-degree extension of $\Q$, and $bN/2$ divides the class number of the real quadratic field $\Q(\sqrt{\delta})$.
\end{thm}

\begin{thm}\label{Thm:Main4}(associated with Theorem \ref{Thm:Murty})
Let $b,N\in \Z$ with $b\ge 1$, $N\ge 3$ and $N\equiv 3 \pmod{4}$. Define
\[f_{N,b}(x)=x^N-x-(N-1)N^b\quad \mbox{and} \quad \delta:=1-N^{(b+1)N-b}.\] Suppose that $f_{N,b}(x)$ is irreducible over $\Q$. Then
\begin{equation}\label{Eq:Mono4}
f_{N,b}(x)\ \mbox{is monogenic if and only if $\delta$ is squarefree.}
\end{equation} Furthermore, no two such monogenic trinomials $f_{N,b}(x)$ generate the same $N$th-degree extension of $\Q$, and the ideal class group of the imaginary quadratic field $\Q(\sqrt{\delta})$ contains an element of order  $(b+1)N-b$.
\end{thm}
\section{Preliminaries}\label{Section:Prelims}
The first theorem is a special case of a result due to Swan \cite{Swan}.
\begin{thm}\label{Thm:Swan}
Let $f(x)=x^N+Ax+B\in \Q[x]$. Then
\[
\Delta(f)=(-1)^{N(N-1)/2}\left(N^{N}B^{N-1}-(-1)^{N}(N-1)^{N-1}A^{N}\right).
\]
\end{thm}

The theorem below is the special case of a theorem due to Jakhar, Khanduja and Sangwan \cite[Theorem 1.1]{JKS2}, which is used to determine the monogenicity of an irreducible trinomial, when applied to trinomials of the form $f(x)=x^N+Ax+B$, where $N\ge 3$.
 We use the notation $q^e\mid \mid n$ to mean that $q^e$ is the exact power of the prime $q$ that divides the integer $n$.

\begin{thm}{\rm \cite{JKS2}}\label{Thm:JKS}
Let $N\ge 3$ be an integer.
Let $K=\Q(\theta)$ be an algebraic number field with $\theta\in \Z_K$, the ring of integers of $K$, having minimal polynomial $f(x)=x^{N}+Ax+B$ over $\Q$. A prime factor $q$ of $\Delta(f)$ does not divide $\left[\Z_K:\Z[\theta]\right]$ if and only if $q$ satisfies one of the following conditions:
\begin{enumerate}[font=\normalfont]
  \item \label{JKS:I1} when $q\mid A$ and $q\mid B$, then $q^2\nmid B$;
  \item \label{JKS:I2} when $q\mid A$ and $q\nmid B$, then
  \[\mbox{either } \quad q\mid A_2 \mbox{ and } q\nmid B_1 \quad \mbox{ or } \quad q\nmid A_2\left(-BA_2^{N}-\left(-B_1\right)^{N}\right),\]
  where $A_2=A/q$ and $B_1=\frac{B+(-B)^{q^j}}{q}$ with $q^j\mid\mid N$;
  \item \label{JKS:I3} when $q\nmid A$ and $q\mid B$, then
  \[\mbox{either } \quad q\mid A_1 \mbox{ and } q\nmid B_2 \quad \mbox{ or } \quad q\nmid A_1\left(-AA_1^{N-1}-\left(-B_2\right)^{N-1}\right),\]
  where $A_1=\frac{A+(-A)^{q^\ell}}{q}$ with $q^\ell\mid\mid (N-1)$, and $B_2=B/q$;
  \item \label{JKS:I5} when $q\nmid AB$, then $q^2\nmid \Delta(f)$.
   \end{enumerate}
\end{thm}

\section{Proof of Theorem \ref{Thm:Main1}}
\begin{proof}
Note that replacing $u$ with $-u$ and $w$ with $-w$ in \eqref{Eq:f} yields the same trinomial $f_{u,w}(x)=x^3-uwx-u^2$. Consequently, we can assume without loss of generality that $u>0$.
 A straightforward calculation shows that
 \[\Delta(f_{u,w})=4u^3w^3-27u^4.\]

 We begin with the proof of \eqref{Eq:Mono1}, and we assume first that $f_{u,w}(x)$ is monogenic. Let $q$ be a prime divisor of $\Delta(f_{u,w})$, and suppose that $f_{u,w}(\theta)=0$. If $u>1$, then condition \eqref{JKS:I1} of Theorem \ref{Thm:JKS} is not satisfied for every prime divisor $q$ of $u$; in which case, $q\mid \left[\Z_K:\Z[\theta]\right]$ and $f(x)$ is not monogenic. Thus, $u=1$ and we can focus on
  \begin{equation}\label{Eq:Newf}
  f_{1,w}(x)=x^3-wx-1\quad \mbox{and} \quad \Delta(f_{1,w})=d=4w^3-27.
  \end{equation}
  Observe then from \eqref{Eq:Newf} that $f_{u,w}(x)$ is irreducible over $\Q$ if and only if $w\not \in \{0,2\}$ by the Rational Root Theorem.

If $q=3$, then $3\mid w$ from \eqref{Eq:Newf}.
         From condition \eqref{JKS:I2} of Theorem \ref{Thm:JKS}, we see that $B_1=0$, which implies that $3\nmid A_2$ since $f(x)$ is monogenic. Hence, we have that $w\equiv \pm 3 \pmod{9}$. However, with $u=1$, examining the two possibilities for $w$ reveals that neither condition \eqref{C2} nor condition \eqref{C3} of Theorem \ref{Thm:KM} can be satisfied. For example, if $w\equiv -3 \pmod{9}$, then $uw\equiv -3\not \equiv 3\pmod{9}$, but $u-w\equiv 4\not \equiv \pm 1 \pmod{9}$.   Thus, $3\nmid \Delta(f_{1,w})$ and $3\nmid w$, so that condition \eqref{C1} of Theorem \ref{Thm:KM} holds. As a consequence, we have shown that no cubic trinomial $f_{1,w}(x)$ in \eqref{Eq:Newf} is monogenic when $3\mid w$.

         Suppose next that $q$ is a prime divisor of $d$ such that $q^2\mid d$. By the previous discussion, we have that $q\ne 3$ and $q\nmid w$. But then condition \eqref{JKS:I5} of Theorem \ref{Thm:JKS} is not satisfied, contradicting the fact that $f(x)$ is monogenic. Hence, $d$ is squarefree, which completes the proof in this direction.

         Conversely, suppose that $u=1$ and $d$ as in \eqref{Eq:Newf} is squarefree. Let $q$ be a prime divisor of $\Delta(f_{1,w})=d$. Then, from \eqref{Eq:Newf}, conditions \eqref{JKS:I1}, \eqref{JKS:I3} and \eqref{JKS:I5} of Theorem \ref{Thm:JKS} are easily seen to be satisfied. If $q\mid w$, then $q=3$ from \eqref{Eq:Newf}, and $3^3\mid d$, contradicting the fact that $d$ is squarefree, which completes the proof of \eqref{Eq:Mono1}.

         Next, let $\FF_1$ be the set as defined in \eqref{FF1}. First note the fact that every element of $\FF_1$ is such that $\Gal_{\Q}(f_{u,w})\simeq S_3$, and that the class number of the unique quadratic subfield $\Q(\sqrt{d})$ of the splitting field of $f_{a,b}(x)$ is divisible by $3$, follows immediately from Theorem \ref{Thm:KM}.

         Let $G(t):=4t^3-27\in \Z[t]$. It is easy to verify that $G(t)$ is irreducible over $\Q$. Then, it follows from a result of Erd\H{o}s \cite{Erdos} that there exist infinitely many integers $z$ such that $G(z)$ is squarefree. (In fact, there exist infinitely many primes $p$ such that $G(p)$ is squarefree \cite{Helfgott}.) Consequently, there exist infinitely many integers $w$ such that $d$ is squarefree, and hence, $\FF_1$ is infinite.

         Finally, to see that no two elements of $\FF_1$ generate the same cubic field, we assume, by way of contradiction, that $f_{1,w_1}(x)\ne f_{1,w_2}(x)$ are elements of $\FF_1$ such that $K_1=\Q(\theta_1)$ and $K_2=\Q(\theta_2)$ are equal, where $f_{1,w_i}(\theta_i)=0$.
          Then, since $f_{1,w_i}(x)$ is monogenic by \eqref{Eq:Mono1}, we have that $\Delta(f_{1,w_i})=\Delta(K_i)$, and so it follows from  \eqref{Eq:Dis-Dis} that
         \begin{equation}\label{Eq:Disequation}
         4w_1^3-27=4w_2^3-27.
         \end{equation} Since $w_1\ne w_2$, it is easy to see that \eqref{Eq:Disequation} has no integer solutions, which completes the proof of the theorem.          \end{proof}

\section{Proof of Theorem \ref{Thm:Main2}}
\begin{proof}
Since $f_{a,b}(x):=x^3-2^{2a}x-3^{b-1}$, we have from Theorem \ref{Thm:Swan} that
\[\Delta(f_{a,b})=\delta=2^{6a+2}-3^{2b+1}.\] We first determine the values of $(a,b)$ for which $f_{a,b}(x)$ is irreducible over $\Q$. If $f_{a,b}(x)$ is reducible over $\Q$, then $f_{a,b}(r)=0$ for some $r\in \Z$, so that $r=\pm 3^c$ for some integer $c\ge 0$ by the Rational Root Theorem. Suppose that $r=-3^c$. Then
\begin{align}
  \nonumber f_{a,b}(-3^c)=0 \Longrightarrow \ & (-3^c)\left((-3)^{2c}-2^{2a}\right)=3^{b-1}\\
 \nonumber\Longrightarrow \ & 3^{2c}-2^{2a}=-3^{b-1-c}\\
 \nonumber\Longrightarrow \ & 3^{2c}+3^{b-1-c}=2^{2a}\\
 \label{Root} \Longrightarrow \ & 3^{b-1-c}(3^{3c-b+1}+1)=2^{2a} \quad \mbox{or} \quad 3^{2c}(1+3^{b-1-3c})=2^{2a}.
\end{align} From the first possibility in \eqref{Root}, we deduce that $c=b-1$. Thus, we get the equation $2^{2a}-3^{2(b-1)}=1$, which has no solutions by Mih\u{a}ilescu's theorem \cite{M}. For the second possibility in \eqref{Root}, we conclude that $c=0$. Hence, this possibility yields the equation $2^{2a}-3^{b-1}=1$, which has only the solution $(a,b)=(1,2)$, again by Mih\u{a}ilescu's theorem \cite{M}. Note then that $r=-1$ and
\[f_{1,2}(x)=(x+1)(x^2-x-3).\] A similar argument shows that no solutions arise from assuming that $r=3^c$. Consequently, $f_{a,b}(x)$ is irreducible over $\Q$ if and only if $(a,b)\ne (1,2)$.

 We now prove \eqref{Eq:Mono2}. With $A=-2^{2a}$, $B=-3^{b-1}$ and $q$ a prime divisor of $\Delta(f_{a,b})=\delta$, we see that only condition \eqref{JKS:I5} is applicable in Theorem \ref{Thm:JKS}. Hence, we conclude that $f_{a,b}(x)$ is monogenic if and only if $q^2\nmid \delta$ for all prime divisors $q$ of $\delta$. That is, $f_{a,b}(x)$ is monogenic if and only if $\delta$ is squarefree, which completes the proof of \eqref{Eq:Mono2}.

Consider next the set $\FF_2$, as defined in \eqref{FF2}. Observe that $\FF_1\cap \FF_2=\varnothing$ since $b\ge 2$. To see that no two elements of $\FF_2$ can generate the same cubic extension, we assume, by way of contradiction, that
\[K_1=\Q(\theta_1)=\Q(\theta_2)=K_2,\]
for some $(a_1,b_1)\ne (a_2,b_2)$, where $f_{a_i,b_i}(\theta_i)=0$. Since each $f_{a_i,b_i}(x)$ is monogenic, then we have from \eqref{Eq:Dis-Dis} that $\delta_1=\delta_2$, where
\[\delta_i:=2^{6a_i+2}-3^{2b_i+1}.\]
Without loss of generality, we can assume that $a_1>a_2$.
Then, by rearranging the equation $\delta_1=\delta_2$, we see that
\[2^{6a_2+2}\left(2^{6(a_1-a_2)}-1\right)=3^{2b_2+1}\left(3^{2(b_1-b_2)}-1\right),\]
which implies that
\begin{equation}\label{Eq:2equations}
3^{2(b_1-b_2)}-2^{6a_2+1}=1 \quad \mbox{and} \quad 2^{6(a_1-a_2)}-3^{2b_2+1}=1.
\end{equation} It is then easy to see that neither equation in \eqref{Eq:2equations} has a solution by Mih\u{a}ilescu's theorem \cite{M}. Hence, no two trinomials in $\FF_2$ generate the same cubic field.

 For each $f_{a,b}(x)\in \FF_2$, we define
\begin{gather}
 \label{uwd} u:=3^{2b-2}, \quad w:=2^{2a}, \quad d:=4uw^3-27u^2=u\delta\\
\nonumber  \mbox{and}\\
  \label{g} g_{a,b}(x):=x^3-uwx-u^2=x^3-2^{2a}3^{2b-2}x-3^{4b-4}.
\end{gather}
Note that
\begin{align*}
3^{3b-3}f_{a,b}(x/3^{b-1})&=3^{3b-3}\left((x/3^{b-1})^3-2^{2a}(x/3^{b-1})-3^{b-1}\right)\\
&=x^3-2^{2a}3^{2b-2}x-3^{4b-4}\\
&=g_{a,b}(x),
\end{align*}
and
\[\Delta(g_{a,b})=3^{6b-6}\left(2^{6a+2}-3^{2b+1}\right)=u^2d=u^3\delta=u^3\Delta(f_{a,b}),\] by Theorem \ref{Thm:Swan}. Moreover, $g_{a,b}(x)$ is irreducible over $\Q$ if and only if $(a,b)\ne (1,2)$, $f_{a,b}(x)$ and $g_{a,b}(x)$ generate the same cubic field, but $g_{a,b}(x)$ is not monogenic since condition \eqref{JKS:I1} of Theorem \ref{Thm:JKS} is not satisfied with the prime $q=3$.
 We see from \eqref{uwd} that $\gcd(u,w)=1$, $d$ is not a square (since $\delta<0$) and $3\nmid w$. Hence, it follows from Theorem \ref{Thm:KM} that $\Gal_{\Q}(f_{a,b})\simeq \Gal_{\Q}(g_{a,b})\simeq S_3$ and that $3\mid h_K$, where $h_K$ is the class number of the imaginary quadratic field
 $K=\Q(\sqrt{d})=\Q(\sqrt{\delta})$. Furthermore, since $\delta=2^{2(3a+1)}-3^{2b+1}<0$ and $(3a+1,2b+1)\ne (2,3)$, it follows from Theorem \ref{Thm:K} that $(2b+1)\mid h_K$. Consequently, we deduce that $(6b+3)/\gcd(3,2b+1)$ divides $h_K$, which completes the proof of the theorem.
\end{proof}
\begin{rem}
We point out that in the proof of Theorem \ref{Thm:Main2} a more elementary approach, such as in \cite{H}, will suffice for every appeal
   to Mih\u{a}ilescu's theorem \cite{M}.
\end{rem}
\section{Proof of Theorem \ref{Thm:Main3}}
\begin{proof}
In Theorem \ref{Thm:Swan}, we have
\begin{equation}\label{Eq:AB Main3}
A=-N(N-1)^{b} \quad \mbox{and}\quad B=-(N-1),
\end{equation} so that
\[\Delta(f_{N,b})=N^N(N-1)^{N-1}((N-1)^{bN}+1)=N^N(N-1)^{N-1}\delta.\]

 We first point out some important facts. Suppose that $q$ is a prime dividing $\gcd(A,\delta)$. Then, $q\mid N$ so that
\[\delta\equiv (-1)^{bN}+1 \equiv 2\pmod{q},\] since $2\mid N$.
But $\delta\equiv 0 \pmod{q}$, and therefore, $q=2$. Since $N\equiv 2 \pmod{4}$, we see that $\delta\equiv 2 \pmod{4}$. Consequently,
\begin{equation}\label{Eq:Fact}
2\mid\mid \delta\quad \mbox{and}\quad \gcd(A,\delta)=2.
\end{equation}

To establish \eqref{Eq:Mono3}, we let $q$ be a prime divisor of $\Delta(f_{N,b})$, and we use \eqref{Eq:AB Main3} and Theorem \ref{Thm:JKS}.
Observe that $q\mid A$ and $q\mid B$ if and only if $q\mid (N-1)$. Hence, condition \eqref{JKS:I1} of Theorem \ref{Thm:JKS} is satisfied for all prime divisors $q$ of $N-1$ if and only if $q^2\nmid (N-1)$, or equivalently, $N-1$ is squarefree.

Next, observe that $q\mid A$ and $q\nmid B$ if and only if $q\mid N$. If $q=2$, then since $N\equiv 2 \pmod{4}$, we have $2\mid \mid N$, and in condition \eqref{JKS:I2} of Theorem \ref{Thm:JKS}, we see that
\begin{gather*}
A_2=\frac{-N(N-1)^{b}}{2} \equiv 1 \pmod{2}\\
\mbox{and}\\
 B_1=\frac{(N-1)^2-(N-1)}{2}=\frac{(N-1)(N-2)}{2}\equiv 0 \pmod{2}.
 \end{gather*} Hence,
 \[A_2(-BA_2^N-(-B_1)^N)\equiv 1 \pmod{2},\]
 so that condition \eqref{JKS:I2} is satisfied for $q=2$.
 Suppose then that $q\mid N$ with $q\ne 2$, and let $N=q^jn$, where $q\nmid n$. Then
 \[A_2=\frac{-N(N-1)^b}{q}\equiv (-1)^{b+1}q^{j-1}n \pmod{q},\]
 \begin{align*}
  B_1&=\frac{(q^jn-1)^{q^j}-(q^jn-1)}{q}\\
  &=\frac{(-1)^{q^j}+q^{2j}n(-1)^{q^j-1}+\sum_{k=2}^{q^j}\binom{q^j}{k}q^{jk}n^k(-1)^{q^j-k}-(q^jn-1)}{q}\\
  &=q^{j-1}n(q^j-1)+\sum_{k=2}^{q^j}\binom{q^j}{k}q^{jk-1}n^k(-1)^{q^j-k}\\
  &\equiv -q^{j-1}n \pmod{q},
\end{align*} and therefore,
\[A_2(-BA_2^N-(-B_1)^N)\equiv (-1)^{b+2}2q^{(j-1)(N+1)}n^{N+1} \not \equiv 0 \pmod{q}\]
if and only if $j=1$. Thus, condition \eqref{JKS:I2} of Theorem \ref{Thm:JKS} is satisfied for all odd prime divisors $q$ of $N$ if and only if $N$ is squarefree.

Since $q\nmid A$ and $q\mid B$ is impossible, we see that condition \eqref{JKS:I3} is not applicable.

If $q\nmid AB$, then $q\ne 2$ and $q\mid \delta$. Hence, condition \eqref{JKS:I5} is satisfied for all prime divisors $q\nmid AB$ if and only if $q^2\nmid \delta$, which is equivalent to $\delta$ being squarefree by \eqref{Eq:Fact}.  Thus, we have established \eqref{Eq:Mono3}.

Suppose that the monogenic trinomials $f_{N,b_1}(x)$ and $f_{N,b_2}(x)$ generate the same $N$th-degree extension of $\Q$. That is, $\Q(\alpha_1)=\Q(\alpha_2)$, where $f_{N,b_1}(\alpha_1)=f_{N,b_2}(\alpha_2)=0$. Then, by \eqref{Eq:Dis-Dis}, we have that
\[(N-1)^{b_1N}+1=(N-1)^{b_2N}+1,\] which implies that $b_1=b_2$ and $f_{N,b_1}(x)=f_{N,b_2}(x)$.

Finally, it follows from Theorem \ref{Thm:AC} that $bN/2$ divides the class number of the real quadratic field $\Q(\sqrt{\delta})$.
\end{proof}
\begin{rem}
  We point out that if $N-1$ is squarefree with $N\ge 6$, then $f_{N,b}(x)=x^N-N(N-1)^bx-(N-1)$ is $p$-Eisenstein for every prime divisor $p$ of $N-1$, and is hence irreducible over $\Q$.
\end{rem}
\section{Proof of Theorem \ref{Thm:Main4}}
\begin{proof}
In Theorem \ref{Thm:Swan}, we have
\begin{equation}\label{Eq:AB}
A=-1 \quad \mbox{and}\quad B=-(N-1)N^b,
\end{equation} so that
\[\Delta(f_{N,b})=(N-1)^{N-1}(1-N^{(b+1)N-b})=(N-1)^{N-1}\delta.\]

To establish \eqref{Eq:Mono4}, let $q$ be a prime divisor of $\Delta(f_{N,b})$. Using \eqref{Eq:AB}, we proceed through Theorem \ref{Thm:JKS} to derive  necessary and sufficient criteria to ensure the monogenicity of $f_{N,b}(x)$. Observe first that $q\mid A$ is impossible so that conditions \eqref{JKS:I1} and \eqref{JKS:I2} of Theorem \ref{Thm:JKS} are not applicable.

If $q\nmid A$ and $q\mid B$, then $q$ must divide $N-1$ since $q\mid \Delta(f_{N,b})$. We see in condition \eqref{JKS:I3} of Theorem \ref{Thm:JKS} that $A_1=0$ and $B_2=-(N-1)N^b/q$. Thus, $q\nmid B_2$ if and only if $q^2\nmid (N-1)$. Hence, condition \eqref{JKS:I3} is satisfied for all prime divisors $q$ of $N-1$ if and only if $N-1$ is squarefree.

If $q\nmid AB$, then $q\mid \delta$. Hence, condition \eqref{JKS:I5} is satisfied for all prime divisors $q\nmid AB$ if and only if $q^2\nmid \delta$, which is equivalent to $\delta$ being squarefree. Since $\delta$ being squarefree implies that $N-1$ is squarefree, we have established \eqref{Eq:Mono4}.

The proof that no two monogenic trinomials $f_{N,b}(x)$ generate the same $N$-th-degree extension of $\Q$ is similar to the corresponding proof in Theorem \ref{Thm:Main3}, and we omit the details.
Finally, it follows from Theorem \ref{Thm:Murty} that the ideal class group of the imaginary quadratic field $\Q(\sqrt{\delta})$ contains an element of order  $(b+1)N-b$.
\end{proof}
\section{Examples to Illustrate the Theorems}
All computations were done using Sage.

\subsection{Theorem \ref{Thm:Main1}}       For $\abs{w}\le 10$, we give in Table \ref{T:1} the values of $w$ and $d$ for trinomials $f_{1,w}(x)$ in $\FF_1$, and the class numbers $h_K$ of the unique quadratic subfield $K=\Q(\sqrt{d})$ of the splitting field of $f_{1,w}(x)$.

\subsection{Theorem \ref{Thm:Main2}}
For the seven smallest squarefree values of $\abs{\delta}$ with $1\le a\le 3$, $2\le b \le 6$ and $\delta<0$, we give in Table \ref{T:2} the values of $a$, $b$, $\delta$, the class number $h$ of the imaginary quadratic field $\Q(\sqrt{\delta})$, and the value of $n:=(6b+3)/\gcd(3,2b+1)$.

\subsection{Theorem \ref{Thm:Main3}}
For the two smallest squarefree values of $\delta$, such that $N(N-1)$ is squarefree, we give in Table \ref{T:3} the values of $N$, $b$, $\delta$, the class number $h$ of the real quadratic field $\Q(\sqrt{\delta})$, and $n:=bN/2$.

\subsection{Theorem \ref{Thm:Main4}}
For the four smallest squarefree values of $\abs{\delta}$, we give in Table \ref{T:4} the values of $N$, $b$, $\delta$, the class number $h$ of the real quadratic field $\Q(\sqrt{\delta})$, and $n:=(b+1)N-b$.

 \begin{table}[h]
 \begin{center}
\begin{tabular}{c|ccccccccccccc}
 $w$ & $-10$ & $-7$ & $-5$ & $-4$ & $-1$ & 1 & 4 & 5 & 7& 8 & 10 \\
  $d$ & $-4027$ & $-1399$ & $-527$ & $-283$ & $-31$ & $-23$ & 229 & 473 &  1345 & 2021 & 3973  \\
  $h_K$ & 9 & 27 & 18 & 3 & 3 & 3 & 3 & 3 & 6 & 3 & 6
  \end{tabular}
\end{center}
\caption{Values of $w$ and $d$ for $f_{1,w}(x)\in \FF_1$ and class numbers of $\Q(\sqrt{d})$}
 \label{T:1}
\end{table}

\begin{table}[h]
 \begin{center}
\begin{tabular}{c|ccccccc}
 $a$ & 1 & 2 & 1 & 2 & 1 & 3 & 2\\
 $b$ & 3 & 4 & 4 & 5 & 5 & 6 & 6\\
 $\delta$ & $-1931$ & $-3299$ & $-19427$ & $-160763$ & $-176891$ & $-545747$ & $-1577939$\\
 $h$ & 21 & 27 & 27 & 66 & 132 & 273 & 624\\
 $n$ & 21 & 9 & 9 & 33 & 33 & 39 & 39
  \end{tabular}
\end{center}
\caption{Values of $a$, $b$, $\delta$, class numbers $h$ of $\Q(\sqrt{\delta})$ and $n$}
 \label{T:2}
\end{table}

\begin{table}[h]
 \begin{center}
\begin{tabular}{ccccc}
 $N$ & $b$ & $\delta$ & $h$ & $n$\\ \hline
 6 & 1 & $15626$ & 24 & 3\\
 6 & 2 & $244140626$ & 1248 & 6
 \end{tabular}
\end{center}
\caption{Values of $N$, $b$, $\delta$, class numbers $h$ of $\Q(\sqrt{\delta})$ and $n$}
 \label{T:3}
\end{table}

 \begin{table}[h]
 \begin{center}
\begin{tabular}{ccccc}
 $N$ & $b$ & $\delta$ & $h$ & $n$\\ \hline
 3 & 2 & $-2186$ & 42 & 7\\
 3 & 3 & $-19682$ & 108 & 9\\
 3 & 4 & $-177146$ & 396 & 11\\
 7 & 1 & $-96889010406$ & 196768 & 13
 \end{tabular}
\end{center}
\caption{Values of $N$, $b$, $\delta$, class numbers $h$ of $\Q(\sqrt{\delta})$ and $n$}
 \label{T:4}
\end{table}

\section*{Acknowledgement} The author thanks the anonymous referee for the helpful suggestions.



\begin{thebibliography}{99}

\bibitem{AC}  N. C. Ankeny and S. Chowla, {\em On the divisibility of the class number of quadratic fields}, Pacific J. Math. {\bf 5} (1955), 321--324.

\bibitem{Cohen} H. Cohen, \emph{A Course in Computational Algebraic Number Theory}, {Springer-Verlag}, 2000.

\bibitem{Erdos} P. Erd\H{o}s, {\em Arithmetical properties of polynomials}, J. London Math. Soc. {\bf 28} (1953), 416--425.

\bibitem{Helfgott}  H. Helfgott, {\em Square-free values of $f(p)$, $f$ cubic}, Acta Math. {\bf 213} (2014), no. 1, 107--135.

\bibitem{H} A. Herschfeld, {\em The equation $2^x-3^y=d$}, Bull. Amer. Math. Soc. {\bf 42} (1936), no. 4, 231--234.

\bibitem{JKS2} A. Jakhar, S. Khanduja and N. Sangwan, \emph{Characterization of primes dividing the index of a trinomial}, Int. J. Number Theory {\bf 13} (2017), no. 10, 2505--2514.

\bibitem{K} Y. Kishi, {\em Note on the divisibility of the class number of certain imaginary quadratic
fields}, Glasgow Math. J. {\bf 51} (2009), 187--191; {\em Corrigendum}, Glasgow Math. J. {\bf 52}
(2010), no. 1, 207--208.

\bibitem{KM} Y. Kishi and K. Miyake, {\em Parametrization of the quadratic fields whose class numbers are divisible by three.} (English summary)
J. Number Theory {\bf 80} (2000), no. 2, 209--217.


\bibitem{M}  P. Mih\u{a}ilescu, {\em Primary cyclotomic units and a proof of Catalan's conjecture}, J. Reine Angew. Math. {\bf 572} (2004), 167--195.

\bibitem{Murty}  M. Ram Murty, {\em The abc conjecture and exponents of class groups of quadratic fields}, Number theory (Tiruchirapalli, 1996), 85--95, Contemp. Math., {\bf 210}, Amer. Math. Soc., Providence, RI, 1998.

\bibitem{Swan}  R. Swan, \emph{Factorization of polynomials over finite fields}, Pacific J. Math. {\bf 12} (1962), 1099--1106.

\end{thebibliography}
\end{document}